\documentclass[12pt]{amsart}
\usepackage{amsfonts,amsthm,amsmath,amssymb}
\usepackage{graphicx}

\newtheorem{thm}{Theorem}[section]
\newtheorem{prop}[thm]{Proposition}
\newtheorem{lemma}[thm]{Lemma}

\newtheorem{remark}[thm]{Remark}

\newtheorem{conj}[thm]{Conjecture}

\def\C{\mathbb{C}}
\def\R{\mathbb{R}}
\def\Z{\mathbb{Z}}
\def\P{\mathbb{P}}

\def\g{\mathfrak{g}}
\def\b{\mathfrak{b}}

\def\t{\mathfrak{t}}

\def\Ocl{\overline{\mathcal O}}

\def\a{\alpha}

\def\phi{\varphi}

\def\G{\Gamma}
\def\l{\lambda}

\def\om{\omega}

\def\Oc{{\mathcal O}}

\def\dim{{\rm dim}}

\renewcommand\o{\overline}

\title{Gelfand-Zetlin polytopes and flag varieties.}
\author{Valentina Kiritchenko}
\date{}\keywords{}


\address{}

\email{valentina.kiritchenko@hcm.uni-bonn.de}

\begin{document}
\maketitle

\bigskip

{\small
Abstract: I construct a correspondence between the Schubert cycles on
the variety of complete flags in $\C^n$ and some faces of the
Gelfand-Zetlin polytope associated with the irreducible representation of $SL_n(\C)$
with a strictly dominant highest weight.  The construction is based on
a geometric presentation of Schubert cells by Bernstein--Gelfand--Gelfand
\cite{BGG} using Demazure modules.
The correspondence between the Schubert cycles and faces is then used to
interpret the classical Chevalley formula in Schubert calculus in
terms of the Gelfand-Zetlin polytopes. The whole picture resembles the picture
for toric varieties and their polytopes.}

\bigskip
\section{Introduction}

Let $G$ be the group $SL_n(\C)$, and $X=G/B$ the flag variety for $G$ (here
$B\subset G$ denotes a Borel subgroup). The main
goal of this paper is to translate to the flag variety some of the rich interplay that
exists between geometry of toric varieties and combinatorics of convex polytopes.
As in the case of toric varieties, there is a polytope $P_H$, namely
{\em Gelfand-Zetlin polytope},
naturally associated with each very ample divisor $H$ on $X$. For a toric variety, an
analogous polytope associated with a divisor $H$ gives information about torus orbits
in the toric variety and their intersection products with $H$. For the flag variety, I will
show how to extract similar information about {\em Schubert cycles} in $X$ and their
intersection products with $H$ using the Gelfand-Zetlin polytope $P_H$.
In particular, the classical Chevalley formula can be reformulated nicely in terms of
Gelfand-Zetlin polytopes (see Theorem \ref{t.Chevalley} and Theorem \ref{t.general}).
Toric and flag varieties are most studied examples of {\em spherical varieties}.
The further goal and motivation for the present paper is to use
the relation between the geometry of flag varieties and Gelfand-Zetlin polytopes
developed here to get new insights into geometry of more general spherical varieties
as outlined in \cite{VK0}.


\medskip

Recall that a {\em Schubert or Bruhat cell} is defined as an orbit of $B$ in $X$ under the
left action, and {\em Schubert cycles} are the cycles in the Chow ring of $X$
represented by the closures of Schubert cells.
Schubert cycles  provide a basis in the Chow ring of $X$, and
the latter is isomorphic to the cohomology ring $H^*(X,\Z)$ of $X$
(see e.g. \cite[1.3]{Brion}).
On the other hand, the cohomology ring of the flag variety is generated by the
degree two classes (see, for instance, \cite[Theorem 3.6.15]{Ma}). The group
$H^2(X,\Z)$ is isomorphic to the Picard group of $X$ and can be identified with
the weight lattice of $G$ so that very ample divisors
are identified with strictly dominant weights (see \cite[1.4.3]{Brion}). Recall that
the weight lattice of $G$ is by definition the character lattice $\Z^{n-1}$ of a
maximal torus in $G$.
The central formula in Schubert calculus is the Chevalley formula
for the intersection product of a Schubert cycle with a divisor (see
Subsection \ref{ss.Chevalley} for more details). The Chevalley formula was proved
independently by Bernstein--Gelfand--Gelfand \cite{BGG} and Demazure
\cite{De} and was already contained in a manuscript of Chevalley \cite{Che}, which for
many years remained unpublished.
This formula allows to express Schubert cycles in terms of divisors thus relating
two different descriptions of the cohomology ring of the flag variety \cite{BGG}.

Fix the upper-triangular Borel subgroup $B^+$. Let $\l$ a strictly dominant
(with respect to $B^+$) weight, and $H_{\l}$ the divisor corresponding
to $\l$. We now assign to $H_\l$ a convex polytope $Q_\l$.
Recall that with each strictly dominant weight $\l$ one can associate the
Gelfand-Zetlin polytope $Q_\l$ (note that Zetlin is sometimes also transliterated as
Cetlin or Tsetlin).  This is a convex polytope in $\R^{d}$ whose vertices
lie in the integral lattice $\Z^d\subset\R^d$
(see Subsection \ref{ss.GZdef} for the definition).
Here $d=n(n-1)/2$  denotes the dimension of $X$. Let $T$ be the diagonal
maximal torus. The integral points inside and at the
boundary of $Q_\l$ parameterize a natural basis of $T$-eigenvectors
(introduced in  \cite{GZ}) in the irreducible representation $V_\l$ of $G$ with the
highest weight $\l$.

I will assign  to each Schubert
cycle in $X$ a face of the Gelfand-Zetlin polytope (see Section
\ref{s.cells2faces}).
My construction depends on a choice of a Borel subgroup $B$ containing the maximal torus $T$
(so in fact, I provide $n!$ different correspondences between Schubert
cycles and faces).  For each choice of $B$,  we first construct a correspondence between
$B$--orbits and faces and then use the one-to-one correspondence between
Schubert cycles and $B$--orbits. The correspondence between $B$-orbits and faces
preserves dimensions. The faces obtained for a given $B$ correspond to
{\em Demazure B-modules} in the representation space $V_\l$.
The freedom in the choice of a Borel subgroup allowed by this
construction is very useful. In many cases, it allows us to choose a face whose
combinatorics captures geometry of a given Schubert cycle especially well (see
Theorem \ref{t.Chevalley} below). It might also lead to an interesting
realization of Schubert calculus in terms of Gelfand-Zetlin polytopes
(this is work in progress with Evgeny Smirnov and Vladlen Timorin). See Section
\ref{s.example} for an example of such calculus in the case $G=SL_3(\C)$.

For a special choice of a Borel subgroup, namely for the lower-triangular Borel
subgroup $B^-$, my construction
gives the correspondence between some of the Schubert cells and faces of the
Gelfand-Zetlin polytope constructed by Kogan using the moment map $X\to Q_\l$
\cite{KThesis} (see Section \ref{s.cells2faces} for more details).
In \cite{KM}, Kogan and Miller extended this correspondence to all Schubert cycles:
they assigned to each Schubert cycle a union of faces
using Caldero's toric degenerations of flag varieties \cite{Caldero}.
Both approaches (with moment map and toric degenerations) only allow to work
with $B^-$--orbits, that is, there is only one way to assign a face or a union of
faces to a given Schubert cycle.

For some of the faces of the Gelfand-Zetlin polytope that correspond to the Schubert
cycles, the Chevalley formula for the intersection product of  a Schubert cycle with
the divisor $H_\l$ admits the following interpretation in terms of the respective face
(cf. Theorem \ref{t.main}). We fix a Borel subgroup $B$ containing $T$, and hence fix a
correspondence between Schubert cycles and faces. Denote by
$\Oc_\G$ the $B$--orbit corresponding to a face $\G$, and by $Z_\G$ the
Schubert cycle defined by $\Oc_\G$. In what follows, we only
consider those faces that do correspond to Schubert cycles. We say that a face $\G$ is
{\em admissible} if for each codimension one orbit $\Oc_\Delta$ in the
closure of the orbit $\Oc_\G$ the face $\G$ contains the face $\Delta$.
In other words, the {\em Bruhat order} on Schubert cycles agrees
with the natural order on faces given by inclusion.
\begin{thm} \label{t.Chevalley} For any admissible face $\G$ we have
$$H_\l Z_\G=\sum_{\Delta\subset\G}d(v,\Delta) Z_\Delta,$$
where the sum is taken over the facets $\Delta$ of $\G$ (that correspond to the
Schubert cells $\Oc_\Delta$ of codimension one at the boundary of $\Oc_\G$).
Here $v$ is a fixed
vertex of the face $\G$ and
$d(v,\Delta)$ denotes the integral distance from  $v$ to the face $\Delta$
(see Section \ref{ss.distance} for the definition).
\end{thm}
Note that in this form the formula is completely analogous to the well-known formula
for toric varieties (e.g. see \cite{VK0}).
There is a generalization of Theorem \ref{t.Chevalley}
that holds for all faces (see Theorem \ref{t.general}).

Many Schubert cycles can be represented by an admissible face
for different choices of $B$, but not all of them. E.g.
for $G=SL_3$ all Schubert cycles can be represented by
admissible faces.  For $G=SL_4$,  exactly two Schubert
cycles can not be represented by an admissible face. These two cycles are given by the
Schubert cells whose closures in the flag variety are not smooth. I conjecture
that all Schubert cycles defined by Schubert cells with smooth closures can be
represented by admissible faces. Note also that if we only take $B^-$
(as in \cite{KThesis,KM}) then already for $SL_3$ there will be a Schubert
cycle such that the corresponding face is not admissible (see Remark \ref{r.two}).

It might be possible to extend the correspondence between Schubert cycles and faces
constructed in this paper to the complete flag varieties for other reductive groups by
replacing the Gelfand-Zetlin polytope with appropriate {\em string polytopes}.

This paper is organized as follows. In Section \ref{s.prelim}, we recall
the definition of the Gelfand-Zetlin polytope and the notion of integral distance.
We also state the classical Chevalley formula. Section \ref{s.cells2faces} contains
the main results: the construction of correspondences
between faces of the Gelfand-Zetlin polytope and Schubert cycles and Chevalley
formula in terms of the Gelfand-Zetlin polytope (Theorem \ref{t.main}). In Section
\ref{s.example}, we consider in detail the example $G=SL_3$. In Section \ref{s.proof},
we study combinatorics of the Gelfand-Zetlin polytope and prove Theorem \ref{t.main}.
We also formulate and prove an extension of Theorem \ref{t.main} to non-admissible
faces (Theorem \ref{t.general}).

\medskip

I am grateful to Michel Brion, Nicolas Perrin and Evgeny Smirnov for useful discussions.
I would also like to thank Jacobs University Bremen, the Hausdorff Center for
Mathematics and the Max Planck Institute for Mathematics in Bonn for hospitality and
support.

\section{Gelfand-Zetlin polytopes and Chevalley formula}\label{s.prelim}
In this section, we recall the definition of the Gelfand-Zetlin polytope and
the Chevalley
formula for the intersection product of a Schubert cycle with a divisor. We also discuss the
notion of integral distance.

\subsection{Gelfand-Zetlin polytope} \label{ss.GZdef}
Let $\l=(\l_1,\ldots,\l_n)$ be a strictly increasing collection of $n$ integer numbers.
To each such collection we assign the irreducible representation
$$\pi_\l:G\to GL(V_\l)$$
with the strictly dominant highest weight
$(\l_2-\l_1)\om_1+\ldots+(\l_n-\l_{n-1})\om_{n-1}$ (which will also be denoted by $\l$),
where $\om_1$,\ldots, $\om_{n-1}$ are the fundamental weights of $G$. To define the
fundamental weights we fix the diagonal maximal torus $T$ and the upper-triangular
Borel subgroup $B^+$.
The Gelfand-Zetlin
polytope  $Q_\l$ associated with $\l$ is a convex polytope in $\R^d$ (recall that
$d=n(n-1)/2$) defined  by the inequalities
$$\begin{array}{ccccccccc}
\l_1&       & \l_2    &         &\l_3     &         &\ldots   &         &\l_n   \\
    &x_{1,1}&         &x_{1,2}  &         & \ldots  &         &x_{1,n-1}&       \\
    &       & x_{2,1} &         &\ldots   &         &x_{2,n-2}&         &       \\
    &       &         & \ddots  &\ldots   &         &         &         &       \\
    &       &         &x_{n-2,1}&         &x_{n-2,2}&         &         &       \\
    &       &         &         &x_{n-1,1}&         &         &         &       \\
\end{array}$$
where $(x_{1,1},\ldots,x_{1,n-1};x_{2,1},\ldots,x_{2,n-2};\ldots;x_{n-2,1},x_{n-2,2};x_{n-1,1})$ are coordinates in $\R^d$ and
 the notation $$\begin{array}{ccc}
                     a &  &b \\
                      & c &
                   \end{array}$$
means $a\le c\le b$. See Figure 1 for a picture of the Gelfand-Zetlin polytope for
$G=SL_3$.

There is a $T$-eigenbasis in $V_\l$ such that its vectors are in one-to-one correspondence
with the integral points inside $Q_\l$ (see for instance \cite[Section 5]{KM} for the description of this
basis). We will denote by the same letter $v$ an integral point in $Q_\l$ and the
corresponding basis vector in $V_\l$.
There is a natural map $p$ that assigns to each integral point $v$ the weight of the
corresponding basis vector $v\in V_\l$. Let us extend this
map by linearity to the map $p:\R^d\to\R^{n-1}$. Denote by $P_\l\subset\R^{n-1}$ the
weight polytope of the representation $V_\l$. The map $p$ sends
the Gelfand-Zetlin polytope $Q_\l$ to the weight polytope $P_\l$
and can be written in coordinates
as follows \cite[2.1.2]{KThesis}.  Let $\a_1$,\ldots, $\a_{n-1}$ be the
simple roots of $G$ (so they form a basis in $\R^{n-1}$ dual with respect to the
Cartan-Killing form to the basis of the fundamental weights $\om_1$,\ldots, $\om_{n-1}$).
Then we have
$$p:(x_{ij})\to (\sum_{i=1}^{n-1} x_{1,i})\a_1+(\sum_{i=1}^{n-2} x_{2,i})\a_2+\ldots+
(x_{n-2,1}+x_{n-2,2})\a_{n-2}+x_{n-1,1}\a_{n-1}+$$
$$+\mbox{ constant vector}.$$

\begin{remark} \label{r.analogous}\em
Note that for any two strictly dominant weights $\l$ and $\mu$ the corresponding
Gelfand-Zetlin polytopes $Q_\l$ and $Q_\mu$ are {\em analogous}, that is, have the same
normal fan. In particular, there is a bijective correspondence between their faces.
This is similar to the toric case, where polytopes corresponding to any two very ample
divisors are analogous.
\end{remark}

\subsection{Integral distance} \label{ss.distance}
Below we recall the notion of {\em integral distance}.
Consider the integral lattice $\Z^d\subset \R^d$ in the affine space $\R^d$.
Let $H$ be a hyperplane  spanned by
lattice vectors, and $v\in\Z^d$ an integral point. Then the
{\em integral distance} $d(v,H)$ from $v$ to the hyperplane $H$  is the index
in $\Z^d$ of the subgroup spanned by the vectors $v-u$ for all $u\in H$.
To compute the integral distance we first find a primitive integral
equation $f(x)=0$ defining $H$, that is, $f(x)=a_0+a_1x_1+\ldots+a_dx_d$ where
$a_i\in\Z$ and the greatest common divisor of $a_0$,\ldots, $a_d$ is $1$. It is then
easy to check that the integral distance between $v$ and $H$ is equal to the absolute
value of $f(v)$.

In the sequel, we will use the notion of integral distance in the following
setting. Let $P$ be a convex lattice polytope of dimension $d$ in $\R^d$.
Recall that a vertex $u$ of $P$ is called {\em simple} if exactly $d$ facets
intersect in $u$ (or equivalently, exactly $d$ edges meet at $u$). In other words,
in the neighborhood of $u$ the polytope $P$ looks like a $d$-dimensional simplex.
Let $\G\subset P$ be a face of $P$, and $\Delta\subset\G$ a facet of $\G$ that
contains at least one simple vertex of $P$. This ensures that
there is a unique hyperplane $H$ such that $H\cap P$ is a facet of $P$ and
$H\cap\G=\Delta$. For any integral point $v\in\G$ we can now define the integral
distance $d(v,\Delta)$ as the integral distance from $v$ to the hyperplane $H$.
Such distances arise naturally in toric geometry when one computes products of
toric orbits with divisors.

\subsection{Bruhat order and Chevalley formula} \label{ss.Chevalley}

Fix a strictly dominant weight $\l$. Recall that $V_\l$ denotes the irreducible
representation with the highest weight $\l$.
We assume that $G/B$ is embedded into the projective space $\P(V_\l)$ as the
$G$--orbit of the line spanned by a highest weight vector $v\in V_\l$. Denote by
$H_\l$ the divisor of hyperplane section on $G/B$ (this is one of the equivalent ways
to identify strictly dominant weights with very ample divisors \cite[1.4]{Brion}).
For each Schubert cell $\Oc$ in $G/B$, the Chevalley formula computes explicitly the
intersection of $H_\l$ with the closure of $\Oc$ as a linear combination of
the closures of Schubert cells at the boundary of $\Oc$.  We will now state this formula.

First, recall
that the choice of a Borel subgroup $B$ in $G$ defines a one-to-one
correspondence between the Schubert cells in $G/B$ and elements of the Weyl group $W$ of
$G$. We identify the Weyl group with
$N(T)/T$, where $N(T)$ is the normalizer of $T$ in $G$.
Then the Schubert cell $\Oc_w$ is the $B$--orbit of the line spanned by
$wv \in V_\l$. Note that the length $l(w)$ (defined as the minimal number of simple
reflections in a decomposition of $w$) is equal to the dimension of $\Oc_w$.
Recall that there is a natural partial order on Schubert cells called {\em Bruhat order}.
We say that $\Oc_{w'}$ {\em precedes} $\Oc_{w}$
with respect to the Bruhat order if $\Oc_{w'}$ is contained in the closure of $\Oc_w$
and $\dim~\Oc_{w'}=\dim~\Oc_w-1$. In other words, $\Oc_{w'}$ is a
boundary divisor in $\Ocl_w$. The Bruhat order can also be defined in terms of the
Weyl group as follows. Denote by $s_\a$ the reflection in the hyperplane perpendicular
to a root $\a$. Then $\Oc_{w'}$ precedes $\Oc_{w}$ if and only if $w'=ws_\a$ for some
root $\a$ and $l(w')=l(w)-1$ (see e.g \cite[Theorem 2.11]{BGG} or
\cite[Proposition 3.6.4]{Ma}).

For each root $\a$, define the linear function $(\cdot,\a)$ (that is, the {\em coroot})
on the weight lattice of $G$ by the property $s_\a\l=\l-(\l,\a)\a$ for all weights $\l$.
(The pairing  $(a,b)$ is often denoted by $\langle a,b^\vee\rangle$
or by $\langle a,b\rangle$.)  Denote
by $Z_w$ the Schubert cycle represented by the closure of the orbit $\Oc_w$.
The following result is
proved in \cite[Proposition 4.1]{BGG} and \cite[Proposition 4.4]{De}:
$$H_\l Z_w=\sum_{\a}(\l,\a)Z_{ws_\a},$$
where the sum is taken over all positive roots $\a$ of $G$ such that
$l(ws_\a)=l(w)-1$. In particular, the coefficients $(\l,\a)$ are always nonnegative.

One of our goals is to interpret this formula in terms of the Gelfand-Zetlin
polytope $Q_\l$. In what follows
we will use the following equivalent
formulation:
$$H_\l Z_w=\sum_{\a}(w\l,\a)Z_{s_\a w},\eqno(2.1)$$
where the sum is taken over all  roots $\a$ such that $w^{-1}\a$ is
positive and  $l(s_\a w)=l(w)-1$.

\section{Correspondence between the Schubert cells and the faces of the Gelfand-Zetlin
polytope.}
\label{s.cells2faces}
In this section, we will construct a correspondence between Schubert cycles and some of
the faces of the Gelfand-Zetlin polytope corresponding to a strictly dominant weight
(by Remark \ref{r.analogous} it does not matter which weight we choose).
\subsection{Schubert cells.} \label{ss.cells} Fix once and for all the diagonal maximal
torus $T\subset G$
and denote by $\t$ its Lie algebra. Everything below (weight vectors, Borel subalgebras
etc.) are assumed to be compatible with $T$. As before we assume that $G/B$ is embedded
into the projectivization $\P(V_\l)$ of the irreducible representation $V_\lambda$
as the $G$--orbit of the line spanned by a highest weight vector.

 We will use the following description of
the Schubert cells from \cite{BGG}.  Let $v\in V_\l$ be a non-zero weight vector
with an {\em extremal weight}.
(Recall that a weight is extremal if it is of the form $w\lambda$ for some
element $w$ in the Weyl group of $G$.) Extremal weights are exactly the vertices of
the weight polytope of $V_\lambda$, and their weight spaces are always one-dimensional.
In what follows, we will not distinguish between non-zero proportional vectors with
the same extremal weight.
Let $B$ be a Borel
subgroup in $G$ containing $T$,  and $\b$ its Lie algebra. Note that all such
Borel subgroups lie in the same orbit under the action of the Weyl group $W$
and there are exactly $|W|$ of them. Denote by $U(\b)$ the universal enveloping
algebra of $\b$. Then the pair $(v,B)$ defines the Schubert cell $\Oc(v,B)$, which is
the $B$-orbit of $v$, and the closure of this cell
in the flag variety can be realized as follows \cite[Lemma 2.12]{BGG}:
$$\o{\Oc(v,B)}=X\cap \P(U(\b)v).$$
Note that $U(\b)v$ is a $B$-invariant vector subspace in $V_\l$ (called
{\em Demazure module}). It would be natural to assign to the cell $\Oc(v,n)$ a face of
the Gelfand-Zetlin polytope $Q_\l$  by taking the convex hull of all basis vectors in the
Gelfand-Zetlin basis that lie in the subspace $U(\b)v$
(we identify the
basis vectors in the Gelfand-Zetlin basis with the integral points in $Q_\l$).
Unfortunately, it might happen that the convex hull is not a face or has wrong dimension.
However, this approach still works after some modification
(see Subsection \ref{ss.faces}).

Two Schubert cells $\Oc(v,B)$ and $\Oc(v',B')$ are conjugate by the action of the
Weyl group (and hence represent the same cohomology class) if and
only if $B'=wB w^{-1}$ and $v'$ is proportional to $wv$ for some element $w\in W$.
If we fix a Borel subgroup $B$ then the Schubert cells $\Oc(wv,B)$ for
all $w\in W$ give the full set of $B$-orbits in the flag variety. In particular,
$$X=\sqcup_{w\in W} \Oc(wv,B).$$
For all possible choices of $v$ and $B$, we get $|W|^2$ Schubert cells forming
$|W|$ orbits under the action of the Weyl group.
Note that there is no canonical identification (i.e. independent of the choice
of a Borel subgroup containing the torus $T$)
between the cohomology classes of the cells $\Oc(v,B)$ and the elements of the
Weyl group. By different choices of a Borel subgroup we assign to
each cohomology class $[\o{ \Oc(v,B)}]$
different elements of the Weyl group that are conjugate to each other.
Since we use simultaneously the cells $\Oc(v,B)$ for
different choices of $B$  we will not identify the Schubert cells with elements of the
Weyl group.

\subsection{Faces of the Gelfand-Zetlin polytope.} \label{ss.faces}
To each cell $\Oc(v,B)$ of dimension $l$ we now assign an $l$-dimensional face
of the Gelfand-Zetlin polytope $Q_\l$.
Recall that to each extremal vector $v$ there corresponds a unique
vertex  of the Gelfand-Zetlin polytope, which we also denote by $v$.
The vertex $v$ is a unique preimage of a vertex of the weight
polytope $P_\l$ under the map $p:Q_\l\to P_\l$. It is easy to show that
the vertex $v$ is simple \cite[2.2.2]{KThesis}.
The edges coming out of $v$ are in one-to-one correspondence with the roots $\a$ of $G$
such that the root space $\g_\a$ does not annihilate the extremal
weight vector $v$ \cite[2.2.3]{KThesis}.
For each such root $\a$ denote by $e(v,\a)$ the edge
corresponding to $\a$. The edge $e(v,\a)$ is uniquely defined by the property that its
projection under the map $p:Q_\l\to P_\l$ is parallel to the root $\a$
(see Section \ref{s.proof} for more details on simple vertices $v$ and edges $e(v,\a)$).

Denote by $R(v,B)$ the set of roots such that the root space
$\g_\a\subset\g$ is contained in $\b$ and does not not annihilate $v$.
The cardinality of $R(v,B)$ is equal to the dimension of the cell $\Oc(v,B)$
\cite[Lemma 2.2]{BGG}.
Let $\{\beta_{i_1},\ldots,\beta_{i_l}\}$ be all roots in the set
$R(v,B)$. Assign to the Schubert cell $\Oc(B,v)$ the $l$-dimensional face
$\G(v,B)$ of the Gelfand-Zetlin polytope spanned by the edges
$e(v,\beta_{i_1}),\ldots,e(v,\beta_{i_l})$.
There is a unique such face since the vertex $v$ is simple. This face can be
thought of as a lifting of the Demazure module $U(\b)v$ to the Gelfand-Zetlin polytope.

\begin{remark} \em
There is an alternative description of the face $\G(v,B)$ using {\em Morse theory}
on polytopes (such analog of the Morse theory was introduced in \cite{Khov}). Namely, choose a
linear function $f_{B}$ on $\R^{n-1}$ that takes positive values on all roots $\a$
whose root spaces $\g_\a$ are contained  in $\b$. Then
the composition $f_{B}\circ p$ is a linear function on the Gelfand-Zetlin polytope.
The face $\G(v,B)$ is then precisely the {\em upper separatrix face} for the function
$f_{B}\circ p$ at the vertex $v$. The {\em upper separatrix face} is by definition the
face spanned by all edges at $v$ going upward with respect to the function
$f_{B}\circ p$ (that is, $f_{B}\circ p$ increases along these edges).

Note that the Bruhat cells $\Oc(v,B)$
can be defined in an analogous way. Namely, there is a Morse function on $X$ given by
the composition of the moment map $X\to P_\l$ with the same function $f_{B}$, and the
cells $\Oc(v,B)$  are the upper separatrix manifolds for this Morse function
(see \cite[Section 4]{Atiyah}).
\end{remark}

We now compare the Bruhat order on the cells $\Oc(v,B)$ with the inclusion order on the faces $\G(v,B)$.
It is easy to see that if $\G(u,B)$ is a facet in $\G(v,B)$, then
$\Oc(u,B)$ precedes $\Oc(v,B)$ with respect to the Bruhat order.
The converse is wrong. I.e. it happens already for $G=SL_3$ that
$\Oc(u,B)$ lies at the boundary of $\Oc(v,B)$ in the flag variety but the face
$\G(u,B)$ does not belong to the face $\G(v,B)$ (see Section \ref{s.example}).
We say that the face $\G(v,B)$ is {\em admissible} if it contains
all faces $\G(u,B)$ such that the Schubert cell $\Oc(u,B)$
precedes the Schubert cell $\Oc(v,B)$ with respect to the Bruhat order.


Denote by $B^-$ the Borel subgroup opposite to the one used to construct the
Gelfand-Zetlin polytope.
If we fix $B=B^-$   and only
vary $v$ then my correspondence between Schubert cycles and faces reduces to the
correspondence defined in \cite{KThesis} (see Remark \ref{r.Kogan}).
Note that the collection of faces assigned to $\Oc(v,B^-)$ in \cite[Section 4]{KM}
always contains $\G(v,B^-)$. In particular, if this collection
consists of just one face (that is, of $\G(v,B^-)$), then the corresponding Schubert
cycle is a Kempf variety \cite[Proposition 2.3.2]{KThesis}, which is a very restrictive
condition (see \cite[Proposition 2.2.1]{KThesis} for a characterization of such
Schubert cycles). Equivalently, the corresponding Schubert polynomial consists of a
single monomial. In particular, it is easy to check that in this case
$\G(v,B^-)$ must be admissible.
An advantage of my construction is that the freedom in the choice of
$B$ allows us to represent more general Schubert cycles by a single admissible
face of the Gelfand-Zetlin polytope (see Remark \ref{r.two}).



 An interesting problem is to describe all admissible faces
as well as the corresponding Schubert cycles.
 It is not true that for each  Schubert cycle  there exists a
 representative $\Oc(v,B)$
such that  the face $\G(v,B)$ is admissible.
There is a counterexample for the flag variety $X_4$ of $SL_4(\C)$. Namely, if the
closure of a Schubert cell $\Oc$ in $X_4$ is not smooth then none of the faces
corresponding to the cohomology class of $\Ocl$ is admissible (there are two such
Schubert cells in $X_4$).
For the flag variety $X_3$ of $SL_3(\C)$ the closure of each Schubert cell is smooth and
every Schubert cycle can be represented by an admissible face (see Section \ref{s.example}).
These examples suggest the
following conjecture:

\begin{conj} If a Schubert cell has smooth closure then its cohomology class can be
represented by an admissible face.
\end{conj}
More generally, suppose that the closure of the Schubert cell $\Oc_w$ corresponding
to an element $w\in W$ has at most $k-1$ irreducible divisors at the boundary
$\Ocl_w\setminus \Oc_w$, where $k$ is the number of pairwise distinct
simple reflections in a reduced decomposition for $w$ (in particular, $k\le n-1$).
This is the case for smooth Schubert cycles by \cite[Proposition 2.2.8]{Brion}).
I conjecture that the Schubert cycle $[\o \Oc_w]$ can be represented by an admissible face.

We now state the Chevalley formula in terms of the
Gelfand-Zetlin polytope. For a weight vector $u$, denote by $p(u)$ the weight of $u$.
\begin{prop} \label{p.distance}
If $\G(u,B)$ is a facet of $\G(v,B)$ (in particular, $p(v)=s_ap(u)$ for some
root $\a$), then
$$|(p(v),\a)|=d(v,\G(u,B)),$$ where $d(v,\G(u,B))$
is the integral distance from $v$ to the
face $\G(u,B)$ as defined in Section \ref{ss.distance}.
\end{prop}
This proposition will be proved in Section \ref{s.proof}.
If we apply it to formula (2.1) we immediately get the following Chevalley formula
for the admissible faces.
\begin{thm}\label{t.main}
If the face $\G(v,B)$ is  admissible  then the Chevalley formula for the Schubert
variety $\o {\Oc(v,B)}$ and the divisor $H_\l$ can be written as
$$H_\l \o {\Oc(v,B)}=\sum d(v,\G(u,B))\o {\Oc(u,B)},$$
where the sum is taken over all Schubert cells $\Oc(u,B)$ that precede $\Oc(v,B)$.
\end{thm}


\section{Example: flag variety for $SL_3(\C)$} \label{s.example}
\begin{figure}
\includegraphics[width=10cm]{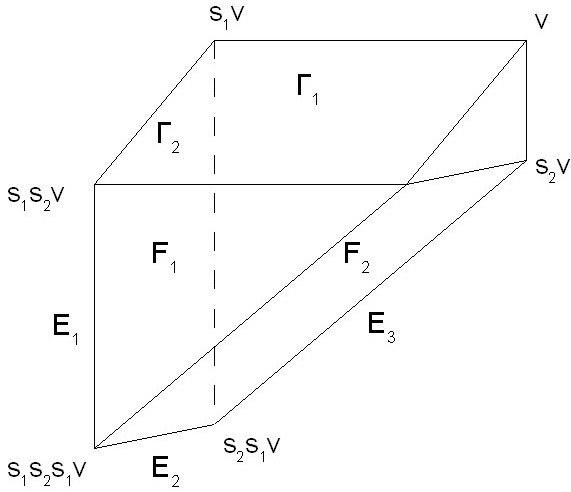}
\caption{}
\end{figure}

Figure 1 shows the Gelfand-Zetlin polytope $Q_\l$ for the irreducible representation of $SL_3(\C)$
with the highest weight
$\l=a\om_1+b\om_2$. This is a polytope in $\R^3$ (with coordinates $x$, $y$ and $z$)
defined by the following six inequalities:
$$0\le x\le a;\quad a\le y\le b; \quad x\le z\le y.$$
The weight polytope $P_\l$ is a hexagon in $\R^2$.
The polytope $Q_\l$ has six simple vertices which are mapped bijectively to the vertices of
the weight polytope $P_\l$ under the map $p$. This bijection is used to label the simple
vertices of $Q_\l$. Namely, we label by $v$ the vertex that goes to the highest weight $\l$. A
simple vertex $u$ is then labeled by $wv$ if $p(u)=wp(v)$ for some element
$w$ from the Weyl group. Put $s_1=s_{\a_1}$
and $s_2=s_{\a_2}$. We denote by $[u_1,u_2]$ the edge of the Gelfand-Zetlin polytope
connecting vertices $u_1$ and $u_2$.

All faces of $Q_\l$ except for a unique non-simple vertex can be represented as
$\G(v,B)$ for some choice of a simple vertex $v$ and a Borel subgroup
$B$. E.g. if $B=B^+$ then $\G(v,B^+)=v$, $\G(s_1v,B^+)$
is the edge $[s_1v,v]$ and
$\G(s_2s_1v,B^+)$ is the face $\{y=b\}\cap Q_\l$. If $B=B^{-}$ then
$\G(v,B^-)=Q_\l$, $\G(s_1v,B^-)$ is the face $\{x=0\}\cap Q_\l$ and
$\G(s_2s_1v,B^-)$ is the edge $[s_2s_1v,s_1s_2s_1v]$.

All faces of $Q_\l$ that do not contain the non-simple vertex are admissible.
In particular, there are two 2-dimensional admissible faces
$\G_1=\G(s_2s_1v,B^+)$ and $\G_2=\G(s_1v,B^-)$ corresponding
to the cells $\Oc(s_2s_1v,B^+)$ and $\Oc(s_1v,B^{-})$.
It is easy to check that these two cells correspond to different
Schubert cycles. Denote these cycles by $Z_{21}$ and $Z_{12}$, respectively (that is,
we label the cohomology class of $\Oc(wv,B^+)$ by $Z_w$ and encode $w=s_1s_2$  by $12$ etc).
There are also six admissible edges that connect simple vertices of $Q_\l$.
These correspond to two Schubert cycles
of dimension one. Namely, the edges $[v,s_1v]$, $[s_1s_2v,s_1s_2s_1v]$ and
$[s_2s_1v,s_2v]$ correspond to $Z_1$, and the other three edges correspond to $Z_2$.
Then Theorem \ref{t.main} applied to the two-dimensional admissible faces tells that
$$H_\l Z_{12}=bZ_1+(a+b)Z_2;\quad H_\l Z_{21}=(a+b)Z_1+aZ_2.$$
\begin{remark}\label{r.two}\em
Note that if we only considered faces $\G(u,B^-)$ for the lower-triangular Borel
subgroup $B^-$ (that is, proceeded as in \cite{KThesis,KM}) then
we would not be able to represent the Schubert cycle $Z_{21}$ by a single admissible
face. Instead, we would get the union of two faces: the rectangular one $\{x=z\}$
and the triangular one $\{y=a\}$. The union of these two faces looks like the admissible
face $\G_1$ (corresponding to $Z_{21}$ by my construction) broken into two pieces.
\end{remark}

We now describe heuristic Schubert calculus on the faces of $Q_\l$. We can represent
Schubert cycle $Z_{21}$ by faces in two different ways: as $\G_1$ and as
$F_1+F_2$,
where $F_1$ and $F_2$ denote the faces given by the equations $y=a$ and $x=z$,
respectively.
The latter representation comes from \cite{KM}. We also represent
$Z_{12}$ by $\G_2$. Finally, we represent the one-dimensional Schubert cycle
$Z_1$ in two ways, by the edge $E_1=[s_1s_2v,s_1s_2s_1v]$ and
the edge $E_3=[s_2v,s_2s_1v]$, and represent $Z_2$ by the edge
$E_2=[s_2s_1v,s_1s_2s_1v]$ (see Figure 1).
We can now compute $Z_{21}Z_{12}$ and $Z_{12}^2$ by intersecting the corresponding
faces:
$$(F_1+F_2)\cap \G_2=E_1+E_2,$$
which is exactly the identity $Z_{21}Z_{12}=Z_1+Z_2$. Similarly,
$$(F_1+F_2)\cap \G_1=E_3$$
gives the identity $Z_{21}^2=Z_1$.
We can also get the identities $Z_1Z_{12}=Z_2Z_{21}=[pt]$ and $Z_1Z_{21}=Z_2Z_{12}=0$
by choosing the edges representing $Z_1$ and $Z_2$ so that they
have transverse intersection with $\G_1$ or $\G_2$. E.g. to find $Z_1Z_{12}$ we
represent $Z_1$ by $E_3$ and $Z_{12}$ by $\G_2$
and get that $\G_2\cap E_3=pt$. Similarly, to find $Z_1Z_{21}$
we  represent $Z_1$ by $E_1$ and $Z_{21}$ by $\G_1$, which yields
$\G_1\cap E_1=\emptyset$.

An analogous Schubert calculus on the Gelfand-Zetlin polytope can be done for arbitrary $n$
\cite{KST}. It can be rigourously justified using the concept of the {\em polytope ring}
whose elements are linear combinations of faces modulo some relations.

\section{Geometry and combinatorics of the Gelfand-Zetlin polytope.} \label{s.proof}
To prove Proposition \ref{p.distance} we have to study the faces of the
Gelfand-Zetlin polytope $Q_\l$. First, we describe explicitly the simple vertices of
$Q_\l$ and the edges going out of simple vertices mostly following \cite{KThesis}. Brief
explanations are provided for the reader's convenience, for more details see \cite[2.1-2.3]{KThesis}.
Next, we will find out under which conditions two
simple vertices are connected by the edge (see Lemma \ref{l.edge}).
Finally, we prove Proposition \ref{p.distance} and formulate and prove a Chevalley formula for
arbitrary faces $\G(v,B)$ (see Theorem \ref{t.general}).

We describe the faces of $Q_\l$ by triangular {\em  diagrams}  following
\cite{KThesis}. Put $x_{0,i}:=\l_i$ for $i=1,\ldots,n$.
It is easy to see that each face of $Q_\l$ is defined by the equations
of the form
$x_{i,j}=x_{i-1,j}$ or $x_{i,j}=x_{i-1,j+1}$ for some $i=1,\ldots,n-1$,
$j=1,\ldots,n-i$.
For a face $\G$, encode all the equations defining $\G$ by the following graph $D(\G)$.
Draw $n$ rows indexed by $1,\ldots,n$ with $n-i+1$ points
$p_{i,1}$,\ldots, $p_{i,n-i+1}$ in the $i$-th row. These are the
vertices of the graph $D(\G)$ (each vertex $p_{i,j}$ corresponds to the coordinate $x_{i-1,j}$).
For each equality $x_{i,j}=x_{i-1,j}$ and $x_{i,j}=x_{i-1,j+1}$ defining
the face $\G$ we draw the edge $e^L_{i+1,j}$ of {\em type $L$}  between the vertices $p_{i+1,j}$ and
$p_{i,j}$ and the edge
$e^R_{i+1,j}$ of {\em type $R$} between   $p_{i+1,j}$ and $p_{i,j+1}$, respectively.
The resulting graph is the {\em diagram} of the face $\G$.
Figure 2 shows the diagrams for the vertices  $v$, $s_1v$  and $s_2v$ of the
Gelfand-Zetlin polytope for $SL_3$ considered in Section \ref{s.example}.

\begin{figure}
\includegraphics[width=10cm]{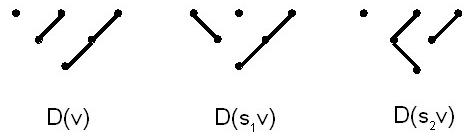}
\caption{}
\end{figure}

\subsection{Simple vertices.} \label{ss.vertices} It is easy to show  that  $v$ is a simple vertex of the
Gelfand-Zetlin polytope if and only if the corresponding diagram $D(v)$ has exactly
$n-i$ edges starting at the $i$-th row and ending  at the $(i+1)$-st row
(for all positive $i<n$) and two such edges never start or end at the same point.
In other words, the graph $D(v)$ is the disjoint union of $n$ simple trees
$T_1(v)$, \ldots,$T_n(v)$.
Each tree $T_i(v)$ starts at the first  row of $D(v)$  and ends at the $i$-th row
(that is, each $T_i$ looks like the Dynkin diagram $A_i$). The vertex of $T_i(v)$ in the first row
will be called the {\em starting point} of $T_i(v)$.
Note that the coordinates $x_{i,j}$ and $x_{k,l}$ of the vertex $v$ are equal if
and only if the vertices $p_{i+1,j}$ and $p_{k+1,l}$ belong to the same tree. The diagram $D(v)$
can also be thought of as an {\em RC-graph} or a {\em pipe dream} (see \cite{KThesis,KM} for details
on connection between pipe-dreams and faces of the Gelfand-Zetlin polytope).
Let us call the diagram of a simple vertex also {\em simple}. There is a different way to characterize
simple diagrams (see \cite[2.2.2]{KThesis}). Namely, the diagram $D(v)$ is simple
if for all $i=2,\ldots,n$ exactly $n-i+1$ edges end at the $i$-th row of $D(v)$, and the edges
$e_{i,j}^L$ are strictly to the left of the edges $e_{i,j}^R$.

Each simple diagram $D(v)$ defines a permutation $\sigma_v$ of elements
$1$,\ldots,$n$ as follows: the vertex $p_{1,i}$ is the starting point of the tree
$T_{\sigma_v(i)}$. It is easy to check that this gives a bijective correspondence between simple
vertices of $Q_\l$ and elements of the symmetric group $S_n$, which is isomorphic to the Weyl group of $G$
(we choose the isomorphism which sends the elementary transposition $(i~(i+1))$ to the simple
reflection $s_{\a_i}$). This bijection is compatible with the bijection between the vertices of
the weight polytope $P_\l$ and elements of the Weyl group, that is, $p(v)=\sigma_v\l$. Indeed,
using the formula for the projection $p:Q_\l\to P_\l$ from Section \ref{ss.GZdef} we get that if
$p(u)=s_{\a_i}p(v)$ (and thus $p(u)=p(v)-(p(v),\a_i)\a_i$), then the sums of coordinates
$\sum_{k=1}^{n-j}x_{j,k}$  for the vertices $v$ and $u$ only differ for $j=i$.
This is only possible if the trees $T_j(v)$ and $T_j(u)$ have the same starting points for all
$j\ne i,(i+1)$, which implies $\sigma_v=s_{\a_i}\sigma_u$.

\subsection{Edges} We now describe the edges of the Gelfand-Zetlin polytope. Let $u$ and $v$ be two
vertices of the Gelfand-Zetlin polytope. We say that the diagram
$D(u)$ is obtained from the diagram $D(v)$ by {\em switching} the edge $e^L_{i,j}$
if the diagrams have the same set of edges with one exception: instead of the edge
$e^L_{i,j}$ the diagram $D(v)$ has the edge $e^R_{i,j}$.
Switching of $e^R_{i,j}$ is defined in the same way. E.g. the diagrams $D(s_1v)$
and $D(s_2v)$ on Figure 2 are obtained from the diagram
$D(v)$ by switching the edges $e^R_{2,1}$ and $e^R_{3,1}$, respectively.
It is easy to see that two vertices $u$ and $v$ are connected by an edge of
the Gelfand-Zetlin polytope only if their diagrams can be obtained from
each other by switching the edge $e^L_{i,j}$ or
$e^R_{i,j}$ for some $i$ and $j$.
If $D(u)$ is obtained from $D(v)$ by
switching the edge $e^L_{i,j}$ then the diagram $D([u,v])$ of the edge
of the Gelfand-Zetlin polytope connecting $u$ and $v$
is obtained from $D(v)$ by deleting the edge $e^L_{i,j}$.

We now focus on the edges going out of a given simple vertex $v$. Their diagrams are
obtained by deleting one of the edges of the diagram $D(v)$.
Denote by $e$ the $i$-th edge of the tree $T_j$ for $i=1,\ldots,j-1$  (that is the
edge of the tree $T_j(v)$ starting at the $i$-th row of the diagram $D(v)$ and
ending at the $(i+1)$-st row). Recall that we denoted by $e(v,\a)$ the edge of
the Gelfand-Zetlin polytope whose
projection $p(e(v,\a))$ is parallel to the root $\a$.
It is easy to check using again the formula for the projection $p:Q_\l\to P_\l$
(see Section \ref{ss.GZdef}) that if we delete the
edge $e$ from the diagram $D(v)$ we get the diagram of the edge $e(v,\a)$, where
$\a=\a_i+\a_{i+1}+\ldots+\a_{j-1}$ if $e$ is of type $L$ and
$\a=-\a_i-\a_{i+1}-\ldots-\a_{j-1}$ if $e$ is of type $R$. Indeed, let $p_{i,s}$ and
$p_{i+1,s}$ be the vertices of the edge $e$. Then switching  $e$
only changes coordinates of $v$ corresponding to the vertices of the tree $T_j(v)$ lying strictly below
$p_{i,s}$. This coordinates increase by the same number $x_{i-1,s+1}(v)-x_{i-1,s}(v)$.
Hence, the sums of coordinates $\sum_{k=1}^{n-r}x_{r,k}$  increase by the same
number for $r=i,\ldots,j-1$, and stay the same for all other $r$. In particular,
for each simple root $\a_i$ the diagram of the edge $e(v,\pm\a_i)$  is obtained from
$D(v)$ by deleting the lowest edge (that is, the $i$-th edge) of the tree $T_{i+1}(v)$,
and the sign in $\pm\a$ is determined by the slope of the lowest edge.
Thus we get an explicit one-to-one correspondence between the edges $e(v,\a)$ of the
Gelfand-Zetlin polytope and the edges of the diagram $D(v)$.

\subsection{Faces $\G(v,B)$ and proof of Proposition \ref{p.distance}} \label{ss.proof}
It is now easy to describe the diagrams of the faces $\G(v,B)$ in terms of
the diagram for $v$. Namely, we should delete all edges in $D(v)$ that correspond to
the roots in $R(v,B)$ under the above correspondence.
E.g. when $B=B^-$ is lower-triangular, the
diagram of $\G(v,B^-)$ is obtained from the diagram $D(v)$ by deleting all edges
of type $R$.

\begin{remark} \label{r.Kogan} The faces $\G(v,B^-)$ are exactly
the so-called {\em Gelfand-Zetlin faces} considered in \cite[Subsection 2.2.1]{KThesis}.
Note that notation in \cite{KThesis} is different: my $x_{i,j}$ is his $\lambda_{i+j,i}$ and my
$\sigma_v$ is his $w^{-1}_v$.
\end{remark}




We now determine under which conditions two simple vertices $u$ and $v$ of the
Gelfand-Zetlin polytope are connected by an edge. The necessary condition
$p(u)=s_\a p(v)$ for some root $\a$ is obviously not sufficient
(e.g. the vertices $s_2v$ and $s_1s_2v$ on Figure 1 are not connected by the
edge though $p(s_1s_2v)=s_1p(s_2v)$).

\begin{lemma} \label{l.edge}
Let $u$ and $v$ be two simple vertices of the Gelfand-Zetlin polytope
such that the weights $p(u)$ and $p(v)$ can be obtained from each other by the reflection
$s_\a$ with respect to some root $\a$. Then $u$ and $v$ are
connected by the edge if and only if the diagram  $D(u)$ can be obtained from
the diagram $D(v)$ by switching the edge of $D(u)$ corresponding to the root
$\a$.
\end{lemma}
\begin{proof} Choose $\a$ so that $(p(v),\a)<(p(u),\a)$.
Then the vertices $v$ and $u$ can only be connected by the edge $e(v,\a)$
(which will then coincide with the edge $e(u,-\a)$), and the lemma
immediately follows from the description of edges in the Gelfand-Zetlin polytope.
\end{proof}


To prove Proposition
\ref{p.distance} we will need the following two lemmas.

\begin{lemma} \label{l.admissible} If $\G(u,B)$ is a facet of $\G(v,B)$,
then the vertices $v$ and $u$ are connected by the edge.
\end{lemma}
\begin{proof} First, note that the assumptions of the lemma imply that $\Oc(u,B)$
precedes $\Oc(v,B)$ with respect to the Bruhat order. Hence, $p(u)=s_\a p(v)$ for some
root $\a\in R(v,B)$. Let $(i~j)$ be the transposition corresponding to $s_\a$, and
$e^R_{i+1,s}$ the edge of the diagram $D(v)$ corresponding to the root
$\a$ (we assume that this edge is of type $R$; type $L$ case is completely analogous).
We now compare the edges starting at the $i$-th rows of the diagrams $D(v)$ and $D(u)$. Let $e^L_{i+1,s-k}$
be the last edge of type $L$ (when going from left to right) starting at the $i$-th row of the diagram
$D(v)$. We want to show that $k=1$, so that $e^R_{i+1,s}$ can be switched and the resulting diagram
remains simple. Consider all edges of $D(v)$ between $e^L_{i+1,s-k}$ and $e^R_{i+1,s}$,
that is, the edges $e^R_{i+1,s-k+1}$,\ldots,$e^R_{i+1,s-1}$.
The above explicit correspondence between simple vertices $v$ and
permutations $\sigma_v$ implies that the trees $T_l(v)$ and $T_l(u)$ have the same starting
points unless $l=i,j$. From this it is easy to deduce that the diagram $D(u)$ contains the edges
$e^L_{i+1,s-k+1}$,\ldots,$e^L_{i+1,s}$. Moreover, if $e^R_{i+1,s-k+l}$ in $D(v)$ for $l=1,\ldots,k-1$
corresponds to a root $\beta$, then
$e^L_{i+1,s-k+l+1}$ in $D(u)$ corresponds to $-\beta$. Finally, $e^L_{i+1,s-k+1}$ corresponds to
the root $-\a$. Hence, the diagrams $D(\G(v,B))$ and  $D(\G(u,B))$
will differ in at least $k$ edges. Indeed, whenever the diagram $D(\G(v,B))$ contains
(or does not contain) the edge $e^R_{i+1,s-k+l}$, the diagram $D(\G(u,B))$ does not contain
(or contains) the edge $e^L_{i+1,s-k+l+1}$ for $l=1,\ldots,k-1$. Also $D(\G(v,B))$ does not
contain the edge $e^L_{i+1,s-k+1}$, while $D(\G(u,B))$ does.
It remains to note that the diagram of $\G(v,B)$ is obtained from the diagram of $\G(u,n)$
by deleting exactly one edge (since $\G(u,B)$ is a facet in $\G(v,n)$). Hence, $k=1$.


\end{proof}
\begin{lemma} \label{l.transp}
If $v$ and $u$ are two simple vertices of the Gelfand-Zetlin polytope such
that $p(v)=s_\a p(u)$ for the root $\a=\a_i+\ldots+\a_{j-1}$, then
$$|(p(v),\a)|=|\l_r-\l_s|,$$
where $s=\sigma^{-1}_v(i)$ and $r=\sigma^{-1}_v(j)$
(that is, $p_{1,s}$ and $p_{1,r}$ are
the starting points of the trees $T_i(v)$ and $T_j(v)$,
respectively).
\end{lemma}
\begin{proof} Since $p(v)=\sigma_v\l$, we have
$(p(v),\a)=(\sigma_v\l,\a)=(\l,\sigma_v^{-1}\a)$.
Note that the reflection defined by the root $\sigma_v^{-1}\a$ corresponds
to the transposition $(\sigma_v^{-1}(i)~\sigma_v^{-1}(j))=(s~r)$.
Hence, $|(\l,\sigma_v^{-1}\a)|=|\l_r-\l_s|$.
\end{proof}
We now prove Proposition \ref{p.distance}. Let $\G(u,B)$ be a facet of $\G(v,B)$.
By Lemma \ref{l.admissible} the vertices $v$ and $u$ are connected
by the edge. We have $p(u)=s_\a p(v)$ for some root $\a$. Suppose that
$\a=\a_i+\a_{i+1}+\ldots+\a_{j-1}$
where $0<i<j<n$. By Lemma \ref{l.transp} we have that $|(p(v),\a)|=|\l_r-\l_s|$,
where $p_{1,s}$ and $p_{1,r}$ are the starting points of the trees $T_i(v)$ and
$T_j(v)$, respectively. We now show that $|\l_r-\l_s|=d(v,\G(u,B))$.
Denote by $e$ the $i$-th
edge of the tree $T_j(v)$.
Since $u$ and $v$ are connected by the edge
we get by Lemma \ref{l.edge} that the diagram $D(u)$ is obtained from $D(v)$
by switching the edge $e$. We again consider the case where $e$ is of type $R$, since
the proof for the other case is completely the same.
Let $p_{i,l+1}$ and $p_{i+1,l}$ be the vertices of the edge $e$. Denote by $F$ the facet
of the Gelfand-Zetlin polytope given by the equation $x_{i-1,l}=x_{i,l}$.
It is easy to check that $\G(u,B)=F\cap\G(v,B)$. Hence,  the integral distance
$d(v,\G(u,B))$ is by definition equal to the distance $d(v,F)$. To compute the latter
we note that the equation $x_{i-1,l}=x_{i,l}$ defining $F$ is already primitive.
Since $p_{i,l}$ belongs to $T_i(v)$ and $p_{i+1,l}$ to $T_j(v)$ we get that the
$x_{i-1,l}$-coordinate of $v$ is equal to $\l_s$ and the $x_{i,l}$-coordinate to $\l_r$.
Hence, $d(v,F)=\l_r-\l_s$.



\subsection{Chevalley formula for arbitrary faces $\G(v,B)$}
The same arguments as in the proof of Proposition \ref{p.distance} allow us to prove
a more general Chevalley type formula for the faces of the Gelfand-Zetlin polytope.
Let $\G(v,B)$ be any (not-necessarily) admissible face, and $\G(u,B)$  a face such that
$\Oc(u,B)$ precedes $\Oc(v,B)$ (but we no longer require that $\G(u,B)\subset\G(v,B)$).
Let $\a=\a_i+\ldots+\a_{j-1}$ be the root such that $p(v)=s_\a p(u)$.
Consider those edges $e_1$,\ldots,$e_k$ ending at the $(i+1)$-st row of the diagram
$D(u)$ that differ by the slope from the corresponding edges at the $(i+1)$-st row of
$D(v)$. Each such edge considered alone gives the diagram of a facet in the
Gelfand-Zetlin polytope. Denote by $F_i$ the facet defined by the edge $e_i$.
Put $d(v,u):=d(v,F_1)+\ldots+d(v,F_k)$.
\begin{thm} \label{t.general} Let $\G(v,B)$ be any (not-necessarily) admissible face. Then
$$H_\l \o {\Oc(v,B)}=\sum d(v,u)\o {\Oc(u,B)},$$
where the sum is taken over all Schubert cells $\Oc(u,B)$ that precede $\Oc(v,B)$.
\end{thm}
Note that for admissible faces
Theorem \ref{t.general} reduces to Theorem \ref{t.main} (since we have $k=1$ by Lemma
\ref{l.admissible} and $\G(u,n)=F_1\cap\G(v,n)$).
Theorem \ref{t.general} is important for a realization of Schubert cycles by unions of
faces of the Gelfand-Zetlin polytope \cite{KST}.
\begin{proof} The proof is almost the same as for admissible faces.
We have $|(p(v),\a)|=|\l_r-\l_s|$ by Lemma \ref{l.transp}. Assume that $r>s$.
We can also write $\l_r-\l_s$ as $(\l_r-\l_{i_{k-1}})+(\l_{i_{k-1}}-\l_{i_{k-2}})+\ldots+(\l_{i_1}-\l_s)$,
where $\l_{i_1}$,\ldots, $\l_{i_{k-1}}$ correspond to the starting points of the trees in
$D(u)$ containing the edges $e_1$,\ldots, $e_{k-1}$, respectively.
It is easy to check that $(\l_{i_{l}}-\l_{i_{l-1}})=d(v,F_l)$ using the same argument as in the proof of
Proposition \ref{p.distance}.
\end{proof}

\end{document}